\documentclass[12pt]{article}

\usepackage[margin=1in]{geometry} 
\usepackage{amsmath} 
\usepackage{amsfonts} 
\usepackage{amsthm} 
\usepackage{mathrsfs}
\usepackage{framed}
\usepackage{graphicx}
\usepackage{graphicx}
\newtheorem{thm}{Theorem}[section]

\begin{document}

\title{The Hockey Stick Theorems in Pascal and Trinomial Triangles}

\author{Sima Mehri
\footnote{Present: PhD student of Mathematics, Sharif University of Technology}
\\Farzanegan High School}

\maketitle

\begin{abstract}

There are some theorems in the Pascal's triangle which their figures resemble to shoot a ball by hockey stick, so they are called hockey stick theorems. P. Hilton and J. Pedersen, in the article "Looking into Pascal Triangle, Combinatorics, Arithmetic and Geometry",\cite{ref2}, have stated the little and big hockey stick and puck theorems in the Pascal's triangle. The big hockey stick theorem is a special case of a general theorem which our goal is to introduce it. We state a hockey stick theorem in the trinomial triangle too.

\end{abstract}

\section{Introduction and Description of Results}
\begin{figure}
\centering
\includegraphics[scale=0.5]{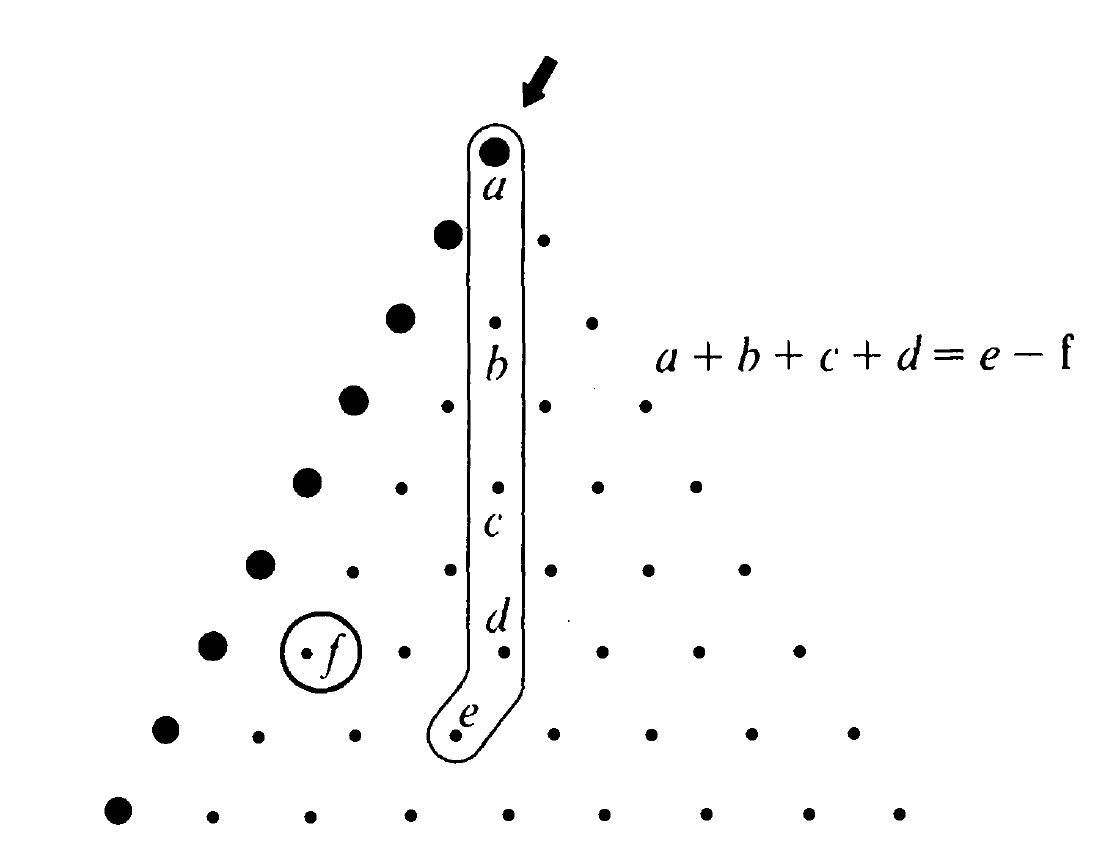}
\caption{The Big Hockey Stick and Puck Theorem}\label{fig:figdozdi}
\end{figure}

The big hockey stick and puck theorem, stated in \cite{ref2} is:
\begin{thm}\cite{ref2}
(The Big Hockey Stick and Puck Theorem)
\[\binom{n}{0}+\binom{n+2}{1}+\binom{n+4}{2}+\binom{n+6}{3}=\binom{n+7}{3}-\binom{n+6}{1}\]
\end{thm}
In \cite{ref2}, this theorem is also demonstrated by Figure \ref{fig:figdozdi}.
We have found the general form of above theorem in Pascal triangle as below.
\begin{thm}\label{pascalhockey}(The Hockey Stick Theorem in Pascal Triangle)
\begin{equation}\label{pascalformula}
\sum _{i=0}^{k}\binom{n+2i}{i} =\sum _{j=0}^{\left\lfloor \frac{k}{2} \right\rfloor }\left(-1\right)^{j}\binom{n+2k-j+1}{k-2j}
\end{equation}
\end{thm}
An example of this theorem is illustrated in Figure \ref{fig:chobechugan dar pascal}.
\setlength{\unitlength}{0.75cm}
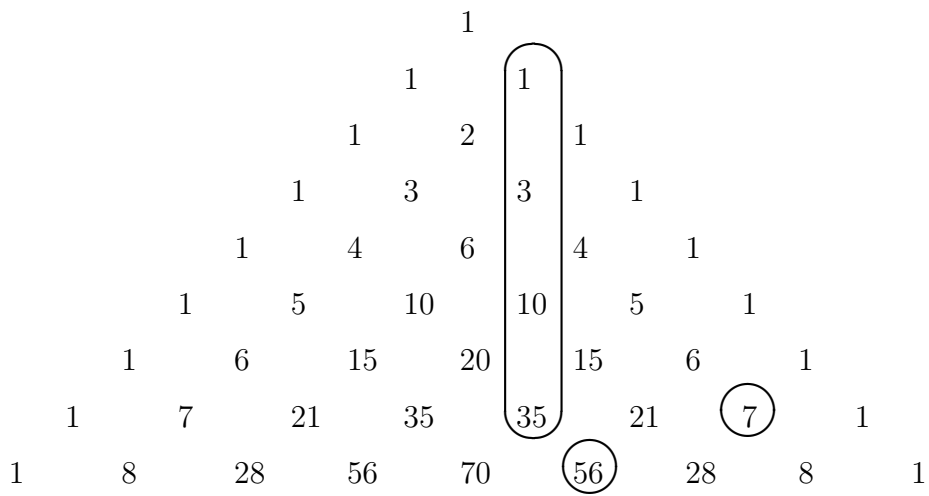
\begin{figure}
\centering
\begin{picture}(16,8)
\linethickness{0.075mm}
\put(8,8){1}
\put(6,6){1}
\put(8,6){2}
\put(10,6){1}
\put(4,4){1}
\put(6,4){4}
\put(8,4){6}
\put(10,4){4}
\put(12,4){1}
\put(2,2){1}
\put(4,2){6}
\put(6,2){15}
\put(8,2){20}
\put(10,2){15}
\put(12,2){6}
\put(14,2){1}
\put(0,0){1}
\put(2,0){8}
\put(4,0){28}
\put(6,0){56}
\thicklines
\put(10.3,0.3){\circle{1}}
\put(8,0){70}
\put(10,0){56}
\put(12,0){28}
\put(14,0){8}
\put(16,0){1}
\put(1,1){1}
\put(3,1){7}
\put(5,1){21}
\put(7,1){35}
\put(9,1){35}
\put(11,1){21}
\put(13,1){7}
\put(13.1,1.3){\circle{1}}
\put(15,1){1}
\put(3,3){1}
\put(5,3){5}
\put(7,3){10}
\put(9,3){10}
\put(11,3){5}
\put(13,3){1}
\put(5,5){1}
\put(7,5){3}
\put(9,5){3}
\put(11,5){1}
\put(7,7){1}
\put(9,7){1}
\thicklines
\put(9.3,4.3){\oval(1,7)}
\end{picture}
\caption{Example of Hocky-Stick:1+3+10+35=56-7}\label{fig:chobechugan dar pascal}
\end{figure}

Now we wish to state the hockey stick theorem in trinomial triangle. First using \cite{ref3} and \cite{ref4}, we explain what is the trinomial triangle.

\begin{figure}
\centering
\begin{picture}(12,6)
\put(0,0){1}
\put(1,0){6}
\put(2,0){21}
\put(3,0){50}
\put(4,0){90}
\thicklines
\put(8.3,0.3){\circle{1}}
\put(5,0){126}
\put(6,0){141}
\put(7,0){126}
\put(8,0){90}
\put(9,0){50}
\put(10,0){21}
\put(10.2,0.3){\circle{1}}
\put(11,0){6}
\put(12,0){1}
\put(12.2,0.3){\circle{1}}
\put(1,1){1}
\put(2,1){5}
\put(3,1){15}
\put(4,1){30}
\put(5,1){45}
\put(6,1){51}
\put(7,1){45}
\put(8,1){30}
\put(9,1){15}
\put(10,1){5}
\put(11,1){1}
\put(2,2){1}
\put(3,2){4}
\put(4,2){10}
\put(5,2){16}
\put(6,2){19}
\put(7,2){16}
\put(8,2){10}
\put(9,2){4}
\put(10,2){1}
\put(3,3){1}
\put(4,3){3}
\put(5,3){6}
\put(6,3){7}
\put(7,3){6}
\thicklines
\put(7.3,3.3){\oval(1,5)}
\put(8,3){3}
\put(9,3){1}
\put(4,4){1}
\put(5,4){2}
\put(6,4){3}
\put(7,4){2}
\put(8,4){1}
\put(5,5){1}
\put(6,5){1}
\put(7,5){1}
\put(6,6){1}
\end{picture}
\caption{Hockey Stick in Trinomial Triangle: $1+2+6+16+45=90-21+1$}\label{fig:trinomial}
\end{figure}
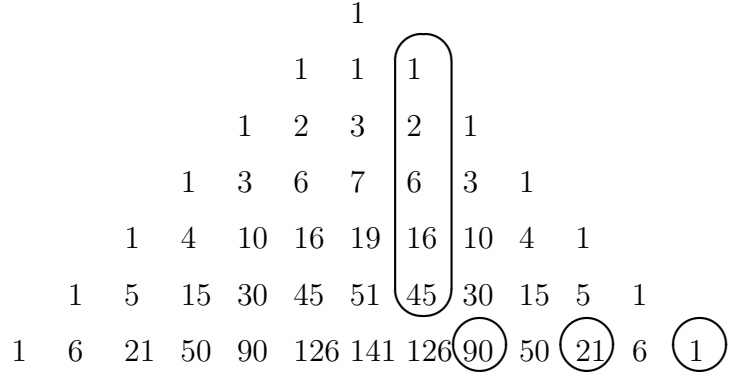

The trinomial triangle is a number triangle of trinomial coefficients. It can be obtained by starting with a row containing a single "1" and the next row containing three 1s and then letting subsequent row elements be computed by summing the elements above to the left, directly above, and above to the right. We show the trinomial triangle in Figure \ref{fig:trinomial}. The trinomial coefficients are placed as Figure \ref{fig:coeff}. Following the notation of Andrews (1990) in \cite{ref1}, the trinomial coefficient  $\binom{n}{k}_2$ with $n\geq 0$ and $-n\leq k\leq n$, is given by the coefficient of $x^{n+k}$ in the expansion of $(1+x+x^2)^n$. Therefore, 

\[\binom{n}{k}_2=\binom{n}{-k}_2\]

Equivalently, the trinomial coefficients are defined by 

\[(1+x+x^{-1})^n=\sum_{k=-n}^{k=n}\binom{n}{k}_2x^k\]
\begin{figure}
\centering
$\begin{array}{ccccccccccc}
& & & & & \binom{0}{0}_2&&&&&\\\\
& & & &\binom{1}{-1}_2&\binom{1}{0}_2&\binom{1}{1}&&&&\\\\
& & & \binom{2}{-2}_2&\binom{2}{-1}_2&\binom{2}{0}_2&\binom{2}{1}&\binom{2}{2}_2&&&\\\\
& & \binom{3}{-3}_2& \binom{3}{-2}_2&\binom{3}{-1}_2&\binom{3}{0}_2&\binom{3}{1}&\binom{3}{2}_2&\binom{3}{3}_2&&\\\\
& \binom{4}{-4}_2& \binom{4}{-3}_2& \binom{4}{-2}_2&\binom{4}{-1}_2&\binom{4}{0}_2&\binom{4}{1}&\binom{4}{2}_2&\binom{4}{3}_2&\binom{4}{4}_2&\\\\
\binom{5}{-5}_2& \binom{5}{-4}_2& \binom{5}{-3}_2& \binom{5}{-2}_2&\binom{5}{-1}_2&\binom{5}{0}_2&\binom{5}{1}&\binom{5}{2}_2&\binom{5}{3}_2&\binom{5}{4}_2&\binom{5}{5}_2\\\\
\end{array}$
\caption{Trinomial Coefficients}\label{fig:coeff}
\end{figure}

We have proven the following theorem in this triangle:
\begin{thm}\label{trinomtheorem}
(The Hockey Stick Theorem in The Trinomial Triangle)
\begin{equation}\label{trinomformula}
\sum _{i=0}^{k}\binom{n+i}{n}_2 =\sum _{s=0}^{\lfloor\frac{k}{2} \rfloor}\left(-1\right)^{s}\binom{n+k+1}{n+2s+1}_2.
\end{equation}
\end{thm} 
For example see Figure \ref{fig:trinomial}.

\section{Proof of Results}
In the proof of both theorems, we use induction.
\begin{proof}
(Theorem \ref{pascalhockey}) We prove this theorem using induction on $k$. By the fact $\binom{n}{n}= \binom{n+1}{n+1}=1$, statement is obvious for the base of induction i.e. $k=1$. Now we assume that the statement for $k$ is true, then the relation \eqref{pascalformula} would be correct. We wish to verify that it is correct for the value $k+1$ too. We have
\begin{align*}
\sum_{i=0}^{k+1}\binom{n+2i}{i} &=\binom{n+2k+2}{k+1} +\sum_{i=0}^{k}\binom{n+2i}{i}=\binom{n+2k+2}{k+1} +\sum _{j=0}^{\left\lfloor \frac{k}{2} \right\rfloor }\left(-1\right)^{j}\binom{n+2k-j+1}{k-2j}\\&=\left[\binom{n+2k+2}{k+1} +\binom{n+2k+1}{k} +\binom{n+2k+1}{k-1}\right]-\\&\quad -\left[\binom{n+2k+1}{k-1} +\binom{n+2k}{k-2} +\binom{n+2k}{k-3}\right]+\\&\quad \quad\vdots \\
&\quad +(-1)^{\left\lfloor \frac{k}{2} \right\rfloor -1}\left[\binom{n+2k-\left\lfloor \frac{k}{2} \right\rfloor +3}{k-2\left\lfloor \frac{k}{2} \right\rfloor +3}+\binom{n+2k-\left\lfloor \frac{k}{2} \right\rfloor +2}{k-2\left\lfloor \frac{k}{2} \right\rfloor +2} +\binom{n+2k-\left\lfloor \frac{k}{2} \right\rfloor +2}{k-2\left\lfloor \frac{k}{2} \right\rfloor +1} \right]+\\&\quad +(-1)^{\left\lfloor \frac{k}{2} \right\rfloor }\left[\binom{n+2k-\left\lfloor \frac{k}{2} \right\rfloor +2}{k-2\left\lfloor \frac{k}{2} \right\rfloor +1} +\binom{n+2k-\left\lfloor \frac{k}{2} \right\rfloor +1}{k-2\left\lfloor \frac{k}{2} \right\rfloor }+\binom{n+2k-\left\lfloor \frac{k}{2} \right\rfloor +1}{k-2\left\lfloor \frac{k}{2} \right\rfloor -1 }\right]+\\&\quad +(-1)^{\left\lfloor \frac{k}{2} \right\rfloor +1}\binom{n+2k-\left\lfloor \frac{k}{2} \right\rfloor +1}{k-2\left\lfloor \frac{k}{2} \right\rfloor -1 }=\\
\intertext{using properties of Pascal triangle, we get}
&=\left[\binom{n+2k+2}{k+1} +\binom{n+2k+2}{k}\right]-\\&\quad -\left[\binom{n+2k+1}{k-1} +\binom{n+2k+1}{k-2}\right]+\\&\quad \quad\vdots \\
&\quad +(-1)^{\left\lfloor \frac{k}{2} \right\rfloor -1}\left[\binom{n+2k-\left\lfloor \frac{k}{2} \right\rfloor +3}{k-2\left\lfloor \frac{k}{2} \right\rfloor +3}+\binom{n+2k-\left\lfloor \frac{k}{2} \right\rfloor +3}{k-2\left\lfloor \frac{k}{2} \right\rfloor +2} \right]+\\&\quad +(-1)^{\left\lfloor \frac{k}{2} \right\rfloor }\left[\binom{n+2k-\left\lfloor \frac{k}{2} \right\rfloor +2}{k-2\left\lfloor \frac{k}{2} \right\rfloor +1} +\binom{n+2k-\left\lfloor \frac{k}{2} \right\rfloor +2}{k-2\left\lfloor \frac{k}{2} \right\rfloor }\right]+\\
&\quad +(-1)^{\left\lfloor \frac{k}{2} \right\rfloor +1}\binom{n+2(k+1)-\left\lfloor \frac{k+1}{2} \right\rfloor +1}{k+1-2\left\lfloor \frac{k+1}{2} \right\rfloor}\mathbf{1}_{\left\lbrace k+1:even\right\rbrace}=\\&=\binom{n+2k+3}{k+1}-\binom{n+2k+2}{k-1}+\cdots +(-1)^{\left\lfloor \frac{k}{2} \right\rfloor }\binom{n+2k-\left\lfloor \frac{k}{2} \right\rfloor +3}{k-2\left\lfloor \frac{k}{2} \right\rfloor +1}+\\&\quad +(-1)^{\left\lfloor \frac{k+1}{2} \right\rfloor }\binom{n+2(k+1)-\left\lfloor \frac{k+1}{2} \right\rfloor +1}{k+1-2\left\lfloor \frac{k+1}{2} \right\rfloor}\mathbf{1}_{\left\lbrace k+1:even\right\rbrace}=\\&=\sum _{j=0}^{\left\lfloor \frac{k+1}{2} \right\rfloor }\left(-1\right)^{j}\binom{n+2(k+1)-j+1}{k+1-2j}
\end{align*}
The statement for $k+1$ is also true, and the proof is completed.
\end{proof}

\begin{proof}(Theorem \ref{trinomtheorem}) To prove this theorem similarly we use induction on $k$. Considering the equation $\binom{n}{n} _2=\binom{n+1}{n+1}_2 $  , the result being immediate if $k=0$. Assuming that the statement for $k$ is true, then the relation \eqref{trinomformula} would be correct. Now we intend to illustrate it is correct for the value $k+1$ too. We have
\begin{align*}
\sum_{i=0}^{k+1}\binom{n+i}{n}_2 &=\binom{n+k+1}{n}_2 +\sum_{i=0}^{k}\binom{n+i}{n}_2=\binom{n+k+1}{n}_2 +\sum_{s=0}^{\left\lfloor \frac{k}{2} \right\rfloor }(-1)^{s} \binom{n+k+1}{n+2s+1} _2=\\&=\left[\binom{n+k+1}{n}_2 +\binom{n+k+1}{n+1}_2 +\binom{n+k+1}{n+2} _2\right]-\\&\quad -\left[\binom{n+k+1}{n+2}_2 +\binom{n+k+1}{n+3}_2 +\binom{n+k+1}{n+4} _2\right]+\\&\quad \quad\vdots \\
&\quad +(-1)^{\left\lfloor \frac{k}{2} \right\rfloor -1}\left[\binom{n+k+1}{n+2\left\lfloor \frac{k}{2} \right\rfloor -2} _2+\binom{n+k+1}{n+2\left\lfloor \frac{k}{2} \right\rfloor -1} _2+\binom{n+k+1}{n+2\left\lfloor \frac{k}{2} \right\rfloor} _2\right]+\\&\quad +(-1)^{\left\lfloor \frac{k}{2} \right\rfloor }\left[\binom{n+k+1}{n+2\left\lfloor \frac{k}{2} \right\rfloor} _2+\binom{n+k+1}{n+2\left\lfloor \frac{k}{2} \right\rfloor +1} _2+\binom{n+k+1}{n+2\left\lfloor \frac{k}{2} \right\rfloor +2} _2\right]+\\&\quad +(-1)^{\left\lfloor \frac{k}{2} \right\rfloor +1}\binom{n+k+1}{n+2\left\lfloor \frac{k}{2} \right\rfloor +2}_2=\\
\intertext{using properties of the trinomial coefficients, we get} &=\sum_{s=0}^{\left\lfloor \frac{k}{2} \right\rfloor }(-1)^{s} \binom{n+k+2}{n+2s+1} _2+(-1)^{\left\lfloor \frac{k}{2} \right\rfloor +1}\binom{n+k+2}{n+2\left\lfloor \frac{k}{2} \right\rfloor +3}_2\\&=\sum_{s=0}^{\left\lfloor \frac{k+1}{2} \right\rfloor }(-1)^{s} \binom{n+k+2}{n+2s+1} _2
\end{align*}
The statement for $k+1$ is also true, and the proof is completed.
\end{proof}
The hockey stick theorem in the trinomial triangles has been proved. This theorem can be translated in Pascal pyramid as follows  :
\[\sum _{i=0}^{k}\; \sum _{2r+s=2n+i} \binom{n+i}{r, s, r-n} =\sum _{j=0}^{\left\lfloor \frac{k}{2} \right\rfloor }\left(\left(-1\right)^{j}\sum _{2r+s=2n+k+2j+2} \binom {n+k+1} {r, s, r-n-2j-1} \right) \] 
Other similar theorems might be obtained for Pascal's four dimensional and even $n$-dimensional pyramid.

\bibliographystyle{plain}

\end{document}